\documentclass{statsoc}
\usepackage{amssymb,amsmath,color,natbib}
\newtheorem{theorem}{Theorem}
\newtheorem{cor}{Corollary}
\newtheorem{prop}{Proposition}
\newtheorem{lemma}{Lemma}

\title[Bayesian change-point problem]{Posterior Convergence and Model Estimation in Bayesian Change-point Problems}
\author[Lian, H.]{Heng Lian}
\address{Division of Mathematical Sciences, SPMS, \\
Nanyang Technological University,\\
Singapore, 637371}
\email{henglian@ntu.edu.sg}
\begin{document}

\begin{abstract}
We study the posterior distribution of the Bayesian multiple change-point regression problem when the number and the locations of the change-points are unknown. While it is relatively easy to apply the general theory to obtain the $O(1/\sqrt{n})$ rate up to some logarithmic factor, showing the exact parametric rate of convergence of the posterior distribution requires additional work and assumptions. Additionally, we demonstrate the asymptotic normality of the segment levels under these assumptions. For inferences on the number of change-points, we show that the Bayesian approach can produce a consistent posterior estimate. Finally, we argue that the point-wise posterior convergence property as demonstrated might have bad finite sample performance in that consistent posterior for model selection necessarily implies the maximal squared risk will be asymptotically larger than the optimal $O(1/\sqrt{n})$ rate. This is the Bayesian version of the same phenomenon that has been noted and studied by other authors.  

\end{abstract}

\section{Introduction}
We consider the regression problem of estimating a piece-wise constant function when the number of segments as well as the locations of its change-points is unknown. This is an old problem that has attracted much  attention recently \citep{goldenshluger06, hariz07, fearnhead08}. Applications of multiple change-point models surged after efficient computations using reversible jump MCMC was discovered \citep{green95}. \citep{green95} applied piece-wise constant function in the study of the coal mining disaster data in the context of Poisson process. A more recent trend of analysis that dispenses with the usage of MCMC for the change-points problem starts with the paper \cite{liulawrence99} where a dynamic programming approach is utilized to marginalize over segment levels and change-point locations. Their original motivation comes from the problem of partitioning DNA sequences into homogeneous segments. This dynamic programming approach is later extended by \cite{fearnhead06, lian08}.

Unlike the above studies, in this paper we are only concerned with the asymptotic properties of Bayesian multiple change-point problems and investigate from the frequentist view the posterior contraction characteristics of a simplified model. Although a piece-wise constant function involves only a finite number of parameters, as we will only consider the case where an upper bound on the number of change-points is available a priori, it is nevertheless best studied from a infinite-dimensional viewpoint and put the estimation problem in the context of function spaces.  Until recently, little is known about the behavior of the posterior distribution of infinite-dimensional models. For consistency issues, \cite{schwartz65} shows that the posterior is consistent when certain tests can be established for the true distribution versus the complement of its neighborhood. \cite{barron99} further developed the theory by sieve construction and metric entropy bounds. Convergence rates are studied in two independent and to some extent overlapping but complementary works (\cite{ghosal00} and \cite{shenwong01}). In particular, \cite{ghosal00} extends the idea of constructing suitable tests in order to bound the convergence rates for both nonparametric and parametric problems, and \cite{ghosal07} further extends the approach to non-i.i.d. observations. The existence of tests for many specific problems can be found in the existing literature although sometimes new tests need to be carefully designed.

In nonparametric Bayesian analysis, we have an i.i.d. sample $Z_1,\ldots,Z_n$ from the distribution $P_0$ with density $p_0$ with respect to some measure on the sample space $(\mathcal{Z},\mathcal{B})$. The model space is denoted by $\mathcal{P}$ which is known to contain the true distribution $P_0$. Given some prior distribution $\Pi$ on $\mathcal{P}$, the posterior is a random measure given by
\[
\Pi^n(A|Z_1,\ldots,Z_n)=\frac{\int_A\prod_{i=1}^np(Z_i)d\Pi(P)}{\int\prod_{i=1}^np(Z_i)d\Pi(P)}\,.
\]
For ease of notation, we will omit the explicit conditioning and write $\Pi^n(A)$ for the posterior distribution. We say that the posterior is consistent if
\[\Pi^n(P\in \mathcal{P} : d(P, P_0) > \epsilon) \rightarrow 0 \mbox{ in } P_0^n \mbox{ probability}\]
for any $\epsilon>0$, where $d$ is some suitable distance function between probability measures.

To study rates of convergence, let $\epsilon_n$ be a sequence decreasing to zero, we say the rate is
at least $\epsilon_n$ if for sufficiently large constant $M$
\[\Pi^n(P : d(P, P_0) > M\epsilon_n)\rightarrow 0 \mbox{ in } P_0^n \mbox{ probability.}\]

We also need a slightly weaker definition of rates of convergence by replacing $M$ with a sequence $M_n$ and requiring that the above posterior mass converge to zero for any sequence $M_n$ that diverges to infinity. This definition is usually required in parametric problems to get rid of the extra $\log n$ factor in the convergence rates. 

In our regression problem, we observe an i.i.d. sample $Z=(Z_1,\ldots,Z_n)$ with the distribution of $Z_i=(X_i, Y_i)$ defined structurally by 
\[ Y_i=\theta_0(X_i)+\epsilon_i\] for  i.i.d Gaussian noise $\epsilon_i\sim N(0,1)$.
and $\theta_0$ is a piece-wise constant function on $[0,1)$ with unknown locations of change-points. We can write $\theta_0(t)=\sum_{j=1}^{k_0} a_j I(t_{j-1}\le t<t_j), t_0=0<t_1<t_2<\ldots<t_{k_0}=1$ using the indicator function and thus $\theta_0$ is parameterized by $(a,t)$, $a=(a_1,\ldots,a_{k_0}), t=(t_0,\ldots,t_{k_0})$.
For simplicity, we assume the marginal distributions for $\{X_i\}$ are i.i.d uniform on $[0,1)$, and note that it is straightforward to extend all the following results to any distribution of $X$ with density bounded away from zero and infinity. Note $P_0$ is fully determined by $\theta_0$ under these assumptions, and thus we also use the space of piece-wise constant functions as our model space which is equivalent to using $\mathcal{P}$. The measure induced by $\theta$ is  denoted by $P_\theta$, and thus $P_{\theta_0}$ is the same as $P_0$, the true distribution.

\cite{lian07b} studied the consistency issue for the above model with the exception that there  $X_i$'s are deterministically chosen on a grid. In that paper, consistency is proved for the case that the true regression function is in the Lipschitz class as well. Another related work is \cite{scricciolo07} where the Bayesian density estimation problem is studied with density approximated by piece-wise constant functions. Besides the fact that they are interested in density estimation instead of regression, the focus of that paper is very different from the current one. They are mostly concerned with the case of approximating a smooth density using step functions and aim to achieve the optimal rates up to a logarithmic factor. For density functions that are piece-wise constant, they prove parametric rate of convergence also with an extra logarithmic factor. A diverging number of grid points is used and thus this approach cannot be used to estimate the number of segments when the density is truly piece-wise constant. 

One simplification of our model compared with those works mentioned at the beginning of this section that focus on the computational issues is that the variance of the noise is assumed to be known here (and actually  $1$ without loss of generality). Investigations of Bayesian regression with unknown noise levels can run into additional technical difficulties especially in the design of appropriate tests. Consistency of a regression problem with unknown noise level is addressed in \cite{choi07}. We hope to be able to address this problem in our context in a future paper. 

In this paper, we focus on the case $\theta_0$ is piece-wise constant and aim to achieve the exact parametric $O(1/\sqrt{n})$ rates of convergence and also study the posterior consistency in the estimation of the number of change-points, which we refer to as the model selection problem. The proofs for the estimation rates involve direct application of general theorems in \cite{ghosal00} but the calculation of the covering number is nontrivial in this case. In order to achieve the exact parametric rate, an additional assumption needs to be made to exclude functions with segment lengths that are too short.

\section{Main results}
Consider the case where we have some a priori bounds for the number of change-points as well as for the segment levels $\{a_j\}$. The model space is defined as
\[\Theta=\{\theta:\theta(t)=\sum_{j=1}^k a_jI(t_{j-1}\le t\le t_j), t_0=0<t_1<\ldots<t_k=1, k\le k_{max}, |a_j|\le K\}\, . \]
By convention, we say $\theta$ with $t_k=1$ has $k$ change-points, which is the same as the number of segments. 

Another equivalent representation of $\Theta$ is $\Theta=\{(a,t)\in [-K,K]^k\times T_k: 1\le k\le k_{max}$, where $T_k$ is the set of $(k+1)$-tuples $(t_0,\ldots,t_k)$ with $t_j<t_{j+1}$. We will not distinguish between these two different representations and $\theta$ can denote either a function or the tuple $(a,t)$. This ambiguity can always be resolved by the context. 

For rates of convergence, the distance $d$ we use is the $L_2$ norm of the function $||\theta||=\left(\int_0^1 \theta^2(x)dx\right)^{1/2}$. Since we only consider uniformly bounded functions, the $L_2$ norm is equivalent to the Hellinger distance (e.g.\cite{ghosal07}, section 7.7). 

We now specify a prior on $\Theta$ using a hierarchical approach. Let $\Theta_k$ be the subspace of $\Theta$ that consists of functions with $k$ change-points, and the prior $\Pi$ is specified as a mixture 
\[\Pi=\sum_{k=1}^{k_{max}} p(k)\Pi_k,\; p(k)>0, \sum_k p(k)=1\]
with $\Pi_k$ the prior measure on $\Theta_k$. We assume that $\Pi_k$ has a density $\pi_k(\theta)$ which can be further decomposed as
\[ \pi_k(a,t)=\pi^a_k(a|t)\pi^t_k(t)\,.\]

The assumption we make on the prior is that 

\textit{(A) The density $\pi^a_k(a|t)$ and $\pi^x_k(t)$ are bounded away from zero and infinity on $[-K,K]^k$ and $T_k$ respectively.}

This assumption is satisfied, for example, when $t_1, t_2,\ldots,t_{k-1}$ are distributed as the order statistic of $k-1$ points uniform distributed on $[0,1)$ while segment levels are independent and uniformly distributed on $[-K,K]$.

The first simple result shows that the posterior rate of convergence is $n^{-1/2}$ up to a logarithmic factor. 
\begin{theorem}\label{lograte} Under assumption (A), the posterior rate of convergence is at least $\epsilon_n=\log^{1/2} n/n^{1/2}$, i.e. $\Pi^n(\theta: ||\theta-\theta_0||>M\epsilon_n)\rightarrow 0$ in $P_0^n$ probability for sufficiently large $M$.
\end{theorem}

Theorem \ref{lograte} considers the convergence rates of the estimation problem. A different problem is the convergence of the posterior for the number of change-points. Under no additional assumptions, we can show that the posterior probability will concentrate on the true number of change-points with probability converging to $1$.

\begin{theorem}\label{modelselection} Under the same assumption (A), we have $\Pi^n(k=k_0)\rightarrow 1$ in $P_0^n$ probability, where $k_0$ is the number of change-points for the true function $\theta_0$. 
\end{theorem}

Nonparametric Bayesian model selection has been investigated in \cite{ghosal08}. The focus of that paper is on conditions under which the adaptive rates are achieved when simultaneously considering models with different rates of contraction. Thus it seems the results presented there cannot be directly applied in our case.

To get rid of the extra logarithmic factor in Theorem \ref{lograte}, one would use the more refined Theorem 2.4 in \cite{ghosal00}, using local covering number instead of the global one. Nevertheless, as shown by Lemma \ref{dbig} in the appendix, the local covering number for $\Theta$ is not bounded as would be required if we set $\epsilon_n=O(1/\sqrt{n})$.
Instead, we consider the smaller model space
\[\Theta^\delta=\{\theta\in\Theta: \min|t_j-t_{j-1}|\ge\delta\}\,.\]
We can define $\Theta_k^{\delta}$ in a similar way and assumption (A) can be modified accordingly. Specification of a prior on $\Theta^\delta$ is easy. Conceptually, we can just restrict $\pi_k^t(t)$ to be supported $\Theta^\delta$ and renormalize the density. Reversible jump algorithms can be easily modified to take into account the constraint. Dynamic programming can also incorporate the pre-specified shortest possible segment length (\cite{lianthesis}. Theorem \ref{lograte} and Theorem \ref{modelselection} is still true on $\Theta^\delta$ with few modifications on the proofs.

\cite{green95} also noticed the practical advantage of avoiding short steps. They proposed using even-numbered order statistics from $2k-1$ uniformly distributed points so that short segment lengths are better penalized. 

As shown in the appendix, putting some lower bound on the segment lengths makes the local covering number bounded by a constant. This requires a very detailed argument to construct the covering. Using this more refined bound on the covering number, we can achieve the exact parametric rate.
\begin{theorem}\label{exactrate}
For any $\delta>0$, under assumption (A), the posterior rate of convergence on $\Theta^\delta$ is at least $\epsilon_n=O(1/\sqrt{n})$. That is, for every $M_n\rightarrow \infty$, we have that $\Pi^n(\theta\in\Theta^\delta: ||\theta-\theta_0||>M_n\epsilon_n)\rightarrow 0$ in $P_0^n$ probability.
\end{theorem}

Combination of Theorem \ref{modelselection} and \ref{exactrate} immediately gives us the rates of convergence for the change-point locations:
\begin{cor}\label{locationrate}
Under the same assumptions as above, the posterior convergence rate for the change-point locations is at least $\epsilon^2_n=O(1/n)$, that is , for any sequence $M_n\rightarrow\infty$, $\Pi^n(\max_{1\le j\le k_0}|t_j-t^0_j|>M_n\epsilon_n^2)\rightarrow 0$ in $P_0^n$ probability.  This rate of course agrees with many frequentist approaches, say using the cumulative sum.
\end{cor}

It is well-known that  the posterior distribution in regular parametric models conditionally converges to a Gaussian distribution under weak conditions. Since the previous results show that the number and locations of the change-points can be consistently estimated, one would naturally conjecture that the posterior distribution for segment levels will converge to a multivariate Gaussian distribution. This is indeed the case as stated in the following theorem:
\begin{theorem}\label{normal}
Suppose the true segment lengths are $l_j=t_j^0-t_{j-1}^0, j=1,\ldots,k_0$. Denoting the posterior distribution of $a=(a_1,\ldots,a_k)$ restricted on the event ${k=k_0}$ (which has a posterior probability converging to $1$) by $\Pi^n_{a|Z}$ and the covariance matrix $I_0={\rm diag}(l_j\cdot n)$, then we have
\[E_{Z|\theta_0}||\Pi^n_{a|Z}-N(\hat{a}(t_0),I_0)||_{TV}\rightarrow 0\, ,\]
where $||P-Q||_{TV}$ is the total variation distance between probability measures $P$ and $Q$, $\hat{a}(t_0)$ is the maximum likelihood estimator for $a$, assuming the true locations of the change-points are known.
\end{theorem}

The above theorems show that the Bayesian procedure possesses very good properties. On the one hand, the exact parametric rate is achieved for the estimation problem in the function space. On the other hand, the number of change-points can be consistently estimated. This is reminiscent of the recent literature on the oracle property of penalized estimators. As shown in \cite{fan01}, the SCAD estimator, when the smoothing parameter is chosen appropriately, can estimate the zero coefficients in a linear regression model as exactly zero with probability converging to one as sample size increases. At the same time, the estimator is still consistent for nonzero coefficients and the asymptotic distribution is the same whether or not the correct zero positions are known. This is called the oracle property by \cite{fan01}. More recently, \cite{leeb08} showed that the oracle property might be misleading in terms of the estimator's finite sample performance and it is impossible to adapt to the unknown zero restrictions without paying a price. The caveat lies in the point-wise nature of the asymptotic theory laid out in \cite{fan01}. The authors of \cite{leeb08} show that an unbounded (normalized) risk results for any estimator possessing the sparsity property. Another related work is \cite{yang05} where the author shows that AIC has a minimax property which cannot be shared with any model selection consistent estimators in a regression problem. 

In our context, similar conclusion can be drawn for the Bayesian multiple change-point problem. Theorem \ref{modelselection} and Theorem {exactrate}  apply to a fixed true piece-wise constant function and thus the convergence as stated is point-wise in nature. It is not difficult to see from the proof of Theorem \ref{exactrate} that the $1/\sqrt{n}$ rate is not actually uniform over the class $\Theta^\delta$. The reason is that to obtain the bound for the local covering number (Lemma \ref{dsmall} in the appendix), the constants involved does depend on $\theta_0$. In particular, the derivation of the lemma requires a lower bound on the size of the jumps of the neighboring segments and thus the convergence is not uniform over $\Theta^\delta$. Intuitively, small jumps makes the estimation more difficult and heavier penalization by the prior must be entertained (possibly by using a prior that depends on the sample size) to achieve model selection consistency at the cost of losing estimation accuracy. As seen in the proof of Theorem \ref{tradeoff}, the difficulty occurs when the size of the jump is of order $O(1/\sqrt{n})$, in which case it becomes difficult to detect the change-point.

Nevertheless, as discussed above, the convergence is uniform if we further restrict our attention on the sub-class:
\[\Theta^{\delta_1,\delta_2}=\{\theta\in\Theta^{\delta_1}, \min |a_j-a_{j-1}|\ge \delta_2\}.\]
We state the uniform convergence as a proposition without proof:
\begin{prop}
For any fixed $\delta_1, \delta_2>0$, the rate of convergence is uniformly at least $\epsilon_n=O(1/\sqrt{n})$. That is, for any $M_n\rightarrow \infty$, $\sup_{\bar{\theta}\in\Theta^{\delta_1,\delta_2}}E_{Z|\bar{\theta}}\Pi^n(\theta\in\Theta^{\delta_1,\delta_2}:||\theta-\bar{\theta}||>M_n\epsilon_n)\rightarrow 0 $.
The property of model selection consistency is still satisfied in this case.
\end{prop}

On the other hand, the following result confirms that we cannot expect the posterior to converge uniformly over the  class $\Theta^{\delta}$ if the method can adapt to the number of change-points. 
Note that the theorem applies for any Bayesian posterior distribution for the change-point problem, not just the specific prior we constructed.
\begin{theorem}\label{tradeoff}
Suppose the  posterior distribution satisfies the model consistency condition: $\Pi^n(k=k_0)\rightarrow 1$ in $P_0^n$ probability, then the maximal $L_2$ convergence of $\theta$ is necessarily slower than the parametric rate $\epsilon_n=O(1/\sqrt{n})$. That is, for some $M_n\rightarrow\infty$,
\[\sup_{\bar{\theta}\in\Theta^{\delta}}E_{Z|\bar{\theta}}\Pi^n(\theta\in\Theta^{\delta}:||\theta-\bar{\theta}||>M_n\epsilon_n)\rightarrow 1 \,.\].
\end{theorem}
The above theorem demonstrated the trade-off between function estimation and model selection for our Bayesian multiple change-point problems. 

\section{Discussion}
In this paper, we investigated in detail some asymptotic properties of Bayesian multiple change-point problems when the noise level is assumed known. We proved estimation rate of convergence as well as model selection consistency of the posterior distribution. 

The main contribution of the paper is to show that the exact parametric rate is achieved for a restricted class of piece-wise constant functions and that this optimal rate cannot be achieved uniformly over the class.

Our theory still leaves some gaps in between. For example, it is still unknown whether it is absolutely necessary to restrict the functions to have not too short segment lengths in order to achieve the optimal rate. The additional restriction makes the local covering number bounded in order to apply the Theorem in \cite{ghosal00}. Besides, the situation with unknown error level is of significant practical importance in which case one should also put a prior on the noise level. The convergence property in this case is still an open problem.

\section*{Appendix}
\subsection*{Some Lemmas}
In preparation for the proofs of the main results, we first collect some lemmas here. Lemma \ref{dbig} below shows that the local covering number is unbounded as remarked in the main text and is not used further in any other proofs. The constant $C$ is used to denote a generic constant which might not be the same at different places. Note that since we are only considering uniformly bounded class of functions, the Hellinger distance, the Kullback-Leibler divergence, as well as the second moment of the likelihood ratio are all equivalent to the $L_2$ norm of the regression function. In the following, we set $\delta_0=\min\{\min_j|t^0_j-t^0_{j-1}|,\min_j|a^0_j-a^0_{j-1}|\}>0$, which bounds the segment lengths as well as the jump size from below.

\begin{lemma} \label{pibig}
Under condition (A), we have the lower bound for the prior concentration when $\epsilon_n\rightarrow 0$,
\[\Pi(\theta: ||\theta-\theta_0||\le \epsilon_n)\ge C p(k_0)\epsilon_n^{3k_0-2}.\]
\end{lemma}
\textit{Proof.} When $\theta=\sum_{j=1}^{k_0} a_jI(t_{j-1}\le t<t_j)\in\Theta_{k_0}$ with $|a_j-a_j^0|<\epsilon_n/2, 1\le j\le k_0$ and $|t_j-t_j^0|<\frac{\epsilon_n^2}{8K^2k_{max}}, 1\le j\le k_0-1$, it is easy to show that $||\theta-\theta_0||^2<\epsilon_n^2$. Since the prior density for $(a,t)$ is bounded away from zero, we get
\[\Pi( \theta: ||\theta-\theta_0||\le \epsilon_n)\ge p(k_0)\Pi_{k_0}(\theta: ||\theta-\theta_0||\le \epsilon_n)\ge Cp(k_0)\epsilon_n^{3k_0-2}.\]

\begin{lemma}\label{pismall}
Let $\delta'=\sqrt{\frac{\delta_0^3}{4k_{max}}}$ ($\delta_0$ is defined immediately before Lemma \ref{pibig}). When $\epsilon<\delta'$, we have that 
$\Pi_k(\theta\in\Theta_k: ||\theta-\theta_0||\le\epsilon)\le C\epsilon^{3k_0-2}, k=1,\ldots,k_{max}$, where $C$ is a constant that depends on $K, k_{max}$ and $\delta_0$. For $k>k_0$, the bound can be refined to $C\epsilon^{3k_0-1}$.
\end{lemma}
\textit{Proof.} First we consider the case $k<k_0$ and $\theta\in\Theta_k$. By the definition of $\delta_0$, the $k_0-1$ intervals $(t_{j}^0-\delta_0/2,t_j^0+\delta_0/2), j=1,\ldots k_0-1$ are nonoverlapping. Thus there is at least one segment of $\theta$ that includes one of these $k_0-1$ intervals. Thus the distance between $\theta$ and $\theta_0$ is at least $\sqrt{\delta(\delta/2)^2}\ge \delta'$, and thus $\Pi_k(\theta\in\Theta_k^\delta: ||\theta-\theta_0||\le\epsilon)=0$. 

When $k\ge k_0, \theta\in\Theta_k$ and $||\theta-\theta_0||<\epsilon$, for any $j$, let $s(j)$ be the index of the interval $[t_{s(j)-1},t_{s(j)})$ which has the largest overlap with $[t_{j-1}^0,t_j^0)$. Obviously the length of the overlap is at least $\delta_0/k_{max}$. This implies $|a_{s(j)}-a_j^0|\le \epsilon\sqrt{\frac{k_{max}}{\delta_0}}$ (otherwise the squared $L_2$ distance between $\theta$ and $\theta_0$ is at least $(a_{s(j)}-a_j^0)^2\frac{\delta_0}{k_{max}}>\epsilon^2$). Similarly, let $t(j)$ be the index of the change-point of $\theta$ that is closest to $t_j^0$, we have $|t_{t(j)}-t_j^0|\le \frac{4\epsilon^2}{\delta_0^2}$ (otherwise the squared distance will be bigger than $\frac{4\epsilon^2}{\delta_0^2}(\frac{\delta_0}{2})^2=\epsilon^2$). The above considerations give us $k_0$ constraints on the segments levels of $\theta$ as well as $k_0-1$ constraints on the change-point locations. Thus under assumption (A), the prior probability $\Pi_k(\theta\in\Theta_k: ||\theta-\theta_0||\le\epsilon)$ is bounded by $C\epsilon^{3k_0-2}$. For refined bound, we consider $k=k_0+1$ only for simplicity. Without loss of generality, we assume $t(j)=j$ and thus $||\theta-\theta_0||\le \epsilon$ implies an additional restriction $(a^0_{k_0}-a_{k_0+1})^2(1-t_{k_0})\le\epsilon^2$. This gives us an additional factor of $\epsilon$ in the bound.

\begin{lemma}\label{dglobal}
$\log D(\epsilon, \Theta)\le b\log(1/\epsilon)+c$, for some constants $b,c>0$ that depends on $K$ and $k_{max}$, where $D(\epsilon, \Theta)$ is the $\epsilon-$covering number of $\Theta$.
\end{lemma}
\textit{Proof.} Choose a grid on the domain $[0,1)$ and another grid on $[-K,K]$
\[ \Delta_t=\left\{\frac{\epsilon^2}{8K^2k_{max}}\cdot i, i\in \mathbb{ N}
\right\}\cap [0,1], \; \Delta_y=\{\epsilon\cdot i, i\in \mathbb{ Z}\}\cap [-K,K]\,.
\]
Let $\tilde{\Theta}=\{\theta\in\Theta, \theta \mbox{ jumps only at points in $\Delta_t$ and takes segment levels in $\Delta_y$}\}$. It is then easy to show that $\tilde{\Theta}$ is an $\epsilon-$covering of $\Theta$ with covering number bounded by 
\[{\lfloor\frac{8K^2k_{max}}{\epsilon^2}\rfloor+1 \choose k_{\max}}(\frac{2K}{\epsilon}+1)^{k_{max}}\,.\]

The next lemma considers the local covering/packing number. In particular, Lemma \ref{dbig} illustrates why we cannot apply Theorem 2.4 from \cite{ghosal00} to obtain the exact parametric rate on $\Theta$.

\begin{lemma}\label{dbig}
\[\log D(\epsilon/2, \{\theta\in\Theta, \epsilon\le ||\theta-\theta_0||\le 2\epsilon)\ge C/\epsilon^2\]
for some constant $C$.
\end{lemma}
\textit{Proof.} Without loss of generality, we consider $\theta_0=0$. We construct a lower bound for the packing number. For simplicity, we assume $1/(4\epsilon^2)$ is an integer. Using the partition of the interval $[0,1)=\cup_{i=1}^{1/4\epsilon^2}[(i-1)4\epsilon^2,i\cdot 4\epsilon^2)$ and construct the piece-wise constant functions $\theta_i(t)=I((i-1)4\epsilon^2\le t<i\cdot 4\epsilon^2)$. Obviously, for this set of functions, we have $||\theta_i||=2\epsilon$ and $||\theta_i-\theta_j||=2\sqrt{2}\epsilon$. The lower bound for the covering number is obtained by the simple relationship between covering number and packing number.

\begin{lemma}\label{dsmall}
For $2\epsilon<\delta'=\sqrt{\frac{\delta_0^3}{4k_{max}}}$,
\[\log D(\epsilon/2, \{\theta\in\Theta^\delta, \epsilon\le ||\theta-\theta_0||\le 2\epsilon)\le C\]
for some constant $C$ that depends on $\delta, \delta_0, K$ and $k_{max}$ but does not depend on $\epsilon$. 
\end{lemma}
\textit{Proof.} Suppose that $||\theta-\theta_0||\le 2\epsilon$. From the proof of Lemma \ref{pismall}, we know that each change-point of $\theta_0$ has a corresponding change-point of $\theta$ that satisfies $|t_{t(j)}-t_j^0|\le 16\epsilon^2/\delta_0^2$. For any segment level $a_j^0$ of $\theta_0$, denote the corresponding index of the segment of $\theta$ that has an overlap of at least $\delta/2$ by $r(j)$, by similar argument as Lemma \ref{pismall}, $|a_j-a_{r(j)}^0|\le 2\sqrt{2}\epsilon/\sqrt{\delta}$. 

To construct a covering, we partition $[0,1)$ into nonoverlapping intervals. In the following, $M,B,N$ are sufficiently large integers to be chosen later. First, each interval $[t_j^0-16\epsilon^2/\delta_0^2,t_j^0+16\epsilon^2/\delta_0^2]$ is partitioned into $M$ subintervals with equal lengths. For the rest of $[0,1)$ we partition it into segments of lengths between $\delta/2B$ and $\delta/B$. Obviously the total number of subintervals does not depend on $\epsilon$. These subintervals falls into two types: (i) the subinterval that contains some change-point of $\theta_0$; (ii) the subinterval that is entirely contained in some segment of $\theta_0$. The function class $F$ that forms a covering is defined as the set of functions which is piece-wise constant with respect to the partition, takes a value of $0$ on type  (i) subintervals and   takes values of the form $a^0_j+\frac{2\sqrt{2}\epsilon}{N\sqrt{\delta}}i,\; i=-N,-(N-1),\ldots, N,$  on type (ii) subintervals if the subinterval is contained in segment $j$ of $\theta_0$.
The size of $F$ is a constant independent of $\epsilon$ and we  show next that it is indeed a $\epsilon/2$-covering.

On subintervals of type (i), the squared $L_2$ distance between $F$ and $\theta$ restricted on these intervals are at most $\frac{16\epsilon^2}{M\delta^2}k_{max}K^2$. Type (ii) subintervals can further be divided into three types: (iii) it contains a change-point of $\theta$ which is closest to some change-point of $\theta_0$; (iv) it contains a change-point of $\theta$ other than those closest to some change-point of $\theta_0$; (v) it is entirely contained in some segment of $\theta$. On subintervals of type (iii) the squared distance  is at most $\frac{16\epsilon^2}{M\delta^2}k_{max}K^2$. On subintervals of type (iv) the squared distance is at most $\frac{\delta}{B}k_{max}(\frac{2\sqrt{2}\epsilon}{\sqrt{\delta}})^2$. On subintervals of type (v) the squared distance is at most $(\frac{2\sqrt{2}\epsilon}{N\sqrt{\delta}})^2$. Thus when $M, B, N$ is large enough, we have a $\epsilon/2$-covering.

\subsection*{Proofs of the main results}

\textit{Proof of Theorem \ref{lograte}.} We apply Theorem 2.1 in \cite{ghosal00} with $\epsilon_n=C\sqrt{\log n/n}$. Condition (2.2) for that theorem is verified in Lemma \ref{dglobal}, condition (2.3) is trivially satisfied and condition (2.4) is verified in Lemmas \ref{pibig}.

\textit{Proof of Theorem \ref{modelselection}.} Theorem \ref{lograte} immediately implies that the under-estimation probablity $\Pi^n(k<k_0)\rightarrow 0$ in $P_0^n$ probability. For over-estimation, it is sufficient to show that $P_0^n(\int_{\Theta_{k_0}}\frac{p_\theta^n(Z)}{p_0^n(Z)}d\pi_{k_0}(\theta)<Cn^{-(3k_0-2+2\xi)/2})\rightarrow 0$ for some $0<\xi<1/2$, and $P_0^n(\int_{\Theta_{k}}\frac{p_\theta^n(Z)}{p_0^n(Z)}d\pi_{k}(\theta)>(\log n)^{-1}n^{-(3k_0-2+2\xi)/2})\rightarrow 0$, when $k>k_0$. 

\textit{step 1.} Let $U_n=\{t\in T_{k_0}:t=t^0+u, u\in R^{k_0+1}, u_0=u_{k_0}=0, |u_i|<c/n\}$ with $\Pi^t_{k_0}(U_n)\ge c'n^{-k_0+1}$, where $\Pi^t_{k_0}$ is the prior measure on the locations of change-points. For any fixed $t\in U_n$, with probability converging to $1$, by considering a small neighborhood of the maximum likelihood estimator $\hat{a}(t)$ for the given $t$ as in Laplace approximation, we have
\[\int\frac{p_\theta^n(Z)}{p_0^n(Z)}d\pi_{k_0}(a|t)
\ge \frac{C}{n^{k_0/2}}\frac{p_{(\hat{a}(t),t)}^n(Z)}{p_0^n(Z)}\ge \frac{C}{n^{k_0/2}}\frac{p_{(a_0,t)}^n(Z)}{p_0^n(Z)}\,.
\]
For any $t\in U_n$, and conditional on $\{X_i\}$, $\log\frac{p_{(a_0,t)}^n(Z)}{p_0^n(Z)}$ is normally distributed with mean $-\frac{1}{2}f(t)^2$ and variance $f(t)^2$, where $f(t)^2=\sum_{j=1}^{k_0-1}(a_{j+1}^0-a_j^0)^2\cdot n_j$, and $n_j$ is the number of $X_i$ that falls into the subinterval $[t_j^0,t_j)$ or $[t_j,t_j^0)$ (depending on the sign of $u_j$). Since $n_j$ is Binomial distributed with mean less than $c$, $f(t)^2=O_p(\xi\log n)$ and thus $\log\frac{p_{(a_0,t)}^n(Z)}{p_0(Z)}\ge -\xi\log n$ with probability converging to 1. Thus, with probability converging to 1, we have $\int_{\Theta_{k_0}}\frac{p_\theta^n(Z)}{p_0^n(Z)}d\pi_{k_0}(\theta)\ge \Pi_{k_0}^t(U_n)\frac{C}{n^{k_0/2}}n^{-\xi}=Cn^{-(3k_0-2+2\xi)/2}$.

\textit{step 2.} Letting $\delta_n={\frac{1}{2\log n n^{3k_0-2(1-\xi)}}}$, and $\epsilon_n=C\log n/\sqrt{n}$, we have that 
\begin{eqnarray*}
&&P_0^n(\int_{\Theta_{k}}\frac{p_\theta^n(Z)}{p_0^n(Z)}d\pi_{k}(\theta)>(\log n)^{-1}n^{-(3k_0-2+2\xi)/2})\\
&\le& P_0^n(\int_{\{||\theta-\theta_0||\le \epsilon_n\}\cap \Theta_k}\frac{p_\theta^n(Z)}{p_0^n(Z)}d\pi_{k}(\theta)>\delta_n)+P_0^n(\int_{\{||\theta-\theta_0||> \epsilon_n\}\cap \Theta_k}\frac{p_\theta^n(Z)}{p_0^n(Z)}d\pi_{k}(\theta)>\delta_n).
\end{eqnarray*}
By the Markov inequality and Fubini's theorem, the first term above is bounded by $\frac{1}{\delta_n}\pi_k(||\theta-\theta_0||\le\epsilon_n)\le \frac{1}{\delta_n}\epsilon_n^{3k_0-1}\rightarrow 0$, where we have made use of Lemma \ref{pismall}.

For the second term, we apply Theorem 1 of \cite{shenwong01} with $\epsilon$ in that theorem replaced by $\epsilon_n$ defined above. Using
\[\int_{\epsilon^2/2^8}^{\sqrt{2}\epsilon}\sqrt{b\log(1/\epsilon)+c}\le\sqrt{b\log\frac{2^8}{\epsilon^2}+c}\cdot \sqrt{2}\epsilon\, ,\]
the entropy condition in that theorem can be verified for $\epsilon=\epsilon_n$. Thus when $C$ is large enough, the second term also converges to $0$.

\textit{Proof of Theorem \ref{exactrate}.}
We apply Theorem 2.4 in \cite{ghosal00} using $\epsilon_n=A/\sqrt{n}$ with $A$ sufficiently large. Condition (2.7) for that theorem is verified in Lemma \ref{dsmall}, and (2.8) is trivially satisfied.  Now we verify (2.9), for which we need to bound $\frac{\Pi(j\epsilon_n\le||\theta-\theta_0||\le 2j\epsilon_n)}{\Pi(||\theta-\theta_0||\le\epsilon_n)}$. When $j<\delta'\sqrt{n}/2A$, $2j\epsilon_n<\delta'$ and Lemma \ref{pismall} can be directly applied to obtain that $\Pi_k(\theta\in\Theta_k^\delta: ||\theta-\theta_0||\le 2j\epsilon_n)\le C(j\epsilon_n)^{3k_0-2}$, and we get $\frac{\Pi(j\epsilon_n\le||\theta-\theta_0||\le 2j\epsilon_n)}{\Pi(||\theta-\theta_0||\le\epsilon_n)}\le Cj^{3k_0-2}\le Cexp(A^2j^2/2)$. For $j\ge \delta'\sqrt{n}/2A$, we bound the probability by $1$, and $\frac{\Pi(j\epsilon_n\le||\theta-\theta_0||\le 2j\epsilon_n)}{\Pi(||\theta-\theta_0||\le\epsilon_n)}\le C(1/\epsilon_n)^{3k_0-2}\le Cexp(A^2j^2/2)$ for this range of $j$.

\textit{Proof of Corollary \ref{locationrate}.} By Theorem \ref{modelselection}, we can assume the number of change-points of $\theta$ is also $k_0$. Then $\max_{1\le i\le k_0}|t_i-t_i^0|>M_n\epsilon_n^2$ implies that $||\theta-\theta_0||^2>(\delta_0/2)^2M_n\epsilon_n^2$. Thus $\Pi_n(\max_{1\le i\le k_0}|t_i-t_i^0|>M_n\epsilon_n^2)\le \Pi_n(\theta\in\Theta^\delta: ||\theta-\theta_0||>(\delta_0/2)\sqrt{M_n}\epsilon_n)\rightarrow 0$.

\textit{Proof of Theorem \ref{normal}.} Fixing one $t\in T_{k_0}^C=\{t\in T_{k_0}: \max_j|t_j-t_j^0|\le C/n\}$, denote the maximum likelihood estimator for $a$ by $\hat{a}(t)$. Let $\Pi^n_{a|t,Z}$ and $\Pi^n_{t|Z}$ be the posterior measure for $a$ conditioning on $t$ and the posterior measure for $t$ respectively. The classical Bernstein-von Mises Theorem implies that $E_0||\Pi^n_{a|t,Z}-N(\hat{a}(t),I_0)||_{TV}\rightarrow 0$.
We have that 
\begin{eqnarray*}
&&E_0||\Pi^n_{a|Z}-N(\hat{a}(t_0),I_0)||_{TV}\\
&\le& ||\int_{T_{k_0}^C} \Pi^n_{a|t,Z}d\Pi^n_{t|Z}-N(\hat{a}(t_0),I_0)||_{TV}+||\int_{(T_{k_0}^C)^c}\Pi^n_{a|t,Z}d\Pi^n_{t|Z}-N(\hat{a}(t_0),I_0)||_{TV}\\
&=&(I)+(II).
\end{eqnarray*}
(I) can be bounded by 
\begin{eqnarray*}
&&E_0||\int_{T_{k_0}^C} \Pi^n_{a|t,Z}d\Pi^n_{t|Z}-N(\hat{a}(t_0),I_0)||_{TV}\\
&\le& E_0\left[\int_{T_{k_0}^C} ||\Pi^n_{a|t,Z}-N(\hat{a}(t),I_0)||_{TV}d\Pi^n_{t|Z}\right]+E_0\left[\int_{T_{k_0}^C} || N(\hat{a}(t),I_0)-N(\hat{a}(t_0),I_0)||_{TV}d\Pi^n_{t|Z}\right].
\end{eqnarray*}
The first term converges to zero by the boundedness of the TV norm and the Fubini's theorem. The second term converges to zero since $||\hat{a}(t)-\hat{a}(t_0)||=o_p(1/\sqrt{n})$. 
Letting $n$ goes to infinity and then $C$ goes to infinity, we see that $E_0||\Pi^n_{a|Z}-N(\hat{a}(t_0),I_0)||_{TV}\rightarrow 0$.

\textit{Proof of Theorem \ref{tradeoff}.} Fix any number $M>0$ and $\gamma>2M$. Define $\theta_0=0$ and $\theta_n=\frac{\gamma}{\sqrt{n}}I(\frac{1}{2}\le t<1)$, a function with a single change-point and jump size $\frac{\gamma}{\sqrt{n}}$. We trivially have $||\theta-\theta_n||\ge \frac{\gamma}{2\sqrt{n}}>\frac{M}{\sqrt{n}}$ for all $\theta\in\Theta_1$ (i.e. $\theta$ is a constant function).  Under $\theta_0$, the posterior probability on $\Theta_1$ converges to $1$ by Theorem \ref{modelselection}. This gives us
\[E_{Z|\theta_0}\Pi^n(\theta: ||\theta-\theta_n||>M/\sqrt{n})\ge 
  E_{Z|\theta_0}\Pi^n(\theta: ||\theta-\theta_n||>M/\sqrt{n}, \theta\in\Theta_1)\ge E_{Z|\theta_0}\Pi^n(\Theta_1)\rightarrow 1.\]
Since the measure $P_0^n$ induces by $\theta_0$ and the measure $P_{\theta_n}^n$ induced by $\theta_n$ are mutually contiguous (this is a straightforward extension of Theorem 7.2 in \cite{vaart98}), we have
\[\sup_{\bar{\theta}\in\Theta^{\delta}}E_{Z|\bar{\theta}}\Pi^n(\theta\in\Theta^{\delta}:||\theta-\bar{\theta}||>M\epsilon_n)\ge E_{Z|\theta_n}\Pi^n(\theta\in\Theta^{\delta}:||\theta-\theta_n||>M\epsilon_n)\rightarrow 1.
\]
Since this is true for any $M$, it is also true for some slowly diverging sequence $M_n$ as in the statement of the theorem.

\bibliographystyle{chicago}
\bibliography{papers.txt,books.txt}

\end{document}